\documentclass[a4paper,12pt,pointlessnumbers,bibtotoc]{article}
\usepackage[a4paper,width=150mm,top=40mm,bottom=50mm,bindingoffset=0mm]{geometry}

\usepackage{blindtext}
\usepackage[UKenglish]{babel}

\usepackage[scaled]{helvet}
\usepackage{lmodern}
\usepackage[T1]{fontenc}
 
\usepackage{graphicx}

\usepackage[babel,german=quotes]{csquotes}
\usepackage[backend=biber ,style=alphabetic,sortcites,natbib=true,sorting=nyt,giveninits=true,backref=true]{biblatex}
\bibliography{bibliography.bib}

\usepackage{url}              
\usepackage{color}            

\usepackage{amssymb, amsmath} 
\usepackage{mathtools}

\usepackage{setspace}         
\newcommand{\fullname}{Emilia Konrad}
\newcommand{\titel}{The Constrained Symplectic Area Functional and its Floer Homology}

\usepackage{amsthm,thmtools,hyperref}
\usepackage[nameinlink,capitalize]{cleveref}    
\hypersetup{
 pdftitle=\titel,
 pdfauthor=\fullname,
 pdfproducer={pdfeTex 3.14159-1.30.6-2.2},
 colorlinks=true,
 pdfborder=0 0 0	
}

\theoremstyle{definition}
\newtheorem{definition}{Definition}[section]
\newtheorem*{notation}{Notation}

\theoremstyle{plain}
\newtheorem{theorem}[definition]{\textup{Theorem}}
\newtheorem{lemma}[definition]{Lemma}
\newtheorem*{remark}{Remark}

\newtheorem{corollary}[definition]{Corollary}

\newtheorem*{claim}{Claim}

\newcounter{assu}
\newtheorem{assumption}[assu]{Assumption}

\newcounter{thmx}
\newtheorem{thm}[thmx]{Theorem}

\crefname{subsection}{Subsection}{Subsections}
\let\refBKP\autoref
\newcommand{\aref}[1]{{\upshape\refBKP{#1}}}

\usepackage{enumitem}
\setlist[enumerate,1]{label={(\roman*)}}

\usepackage{tikz-cd}
\usepackage{caption}
\usepackage{subcaption}
\usepackage{wrapfig}
\usepackage{graphicx, array, blindtext}

\newcommand{\nach}{\longrightarrow}             
\newcommand{\R}{\mathbb{R}}                     
\newcommand{\N}{\mathbb{N}}                     
\newcommand{\inv}{^{-1}}                        
\newcommand{\rund}[1]{\left(#1\right)}          
\renewcommand{\d}{\operatorname{d}\!}           
\newcommand{\norm}[1]{\left\|#1\right\|}        
\renewcommand{\epsilon}{\varepsilon}            
\newcommand{\menge}[1]{\left\{#1\right\}}       
\renewcommand{\phi}{\varphi}                    
\renewcommand{\l}{\ell}                         
\newcommand{\tang}[2]{T_{#1}#2}                 
\newcommand{\bet}[1]{\left|#1\right|}           
\newcommand{\auf}{\longmapsto}                  

\usepackage{mathrsfs}                           
\renewcommand{\L}{\mathscr{L}}                  
\newcommand*{\loops}[2][]{\L^{#1}#2}            
\renewcommand{\S}{{\mathbb{S}^1}}                 
\renewcommand{\star}{^*}                        
\renewcommand{\H}{\mathcal{H}}                  
\newcommand{\intlim}[1]{\int\limits_{#1}}       
\newcommand{\folge}[2]{\rund{#1_{#2}}_{#2\in\N}}    

\newcommand{\liou}{\Lambda}                     
\renewcommand{\vec}[1]{\Gamma T#1}              
\newcommand{\act}[1]{\mathfrak{a}^{#1}}         
\newcommand{\Act}[1]{\mathcal{A}^{#1}}          
\newcommand{\crit}[1]{\operatorname{crit}(#1)}  
\newcommand{\M}{\mathcal{M}}                    
\newcommand{\dt}{\d\,}                          
\newcommand{\dbei}[2]{\left.\frac{\d}{\d#1}\right|_{#1=#2}} 
\newcommand{\ddt}[1]{\frac{\d}{\d#1}}           
\newcommand{\pbei}[2]{\left.\frac{\partial}{\partial#1}\right|_{#1=#2}} 

\newcommand{\C}{\mathbb{C}}                     

\newcommand{\F}{\mathcal{F}}                    
\newcommand{\W}{\mathcal{W}}                    
\newcommand{\ind}{\operatorname{ind}}           
\newcommand{\eps}{\epsilon}                     
\newcommand{\rs}{\R\times\S}                    
\newcommand{\G}{\mathcal{G}}                    

\let\oldforall\forall
\renewcommand{\forall}{\,\oldforall}            
\let\oldexists\exists                          
\renewcommand{\exists}{\,\oldexists}            

\newcommand{\lag}{\tau}                         
\newcommand{\mult}{\chi}                        

\newcommand{\inach}{\overset{\sim}{\nach}}      
\newcommand{\id}{\operatorname{id}}             
\let\oldqed\qed
\renewcommand{\qed}{\oldqed\bigskip}            
\newcommand{\coker}{\operatorname{coker}}       

\newcommand{\ndo}[1]{\operatorname{End}\rund{#1}}   
\newcommand{\im}{\operatorname{im}}             
\newcommand{\Z}{\mathbb{Z}}                     

\newcommand{\hess}{\operatorname{Hess}}         

\newcommand{\spec}{\operatorname{spec}}         

\newcommand{\X}{\mathcal{X}}                    
\newcommand{\Y}{\mathcal{Y}}                    

\newcommand{\cnt}[1]{\bet{#1}_2}                

\usepackage{bbm}    

\newcommand{\bott}{\mu}                             

\newcommand{\Mb}{\widetilde{\M}}                    

\begin{document} 
\pagenumbering{roman}
\begin{center}
    \begin{large}
        \vspace*{-2cm}\textbf{The Constrained Symplectic Area Functional\\ and its Floer Homology}
    \end{large}
      
    \vspace{0.4\baselineskip}
    \fullname
     \vspace{-0.2\baselineskip}
    \begin{abstract} This paper introduces a new Floer homology for periodic Reeb orbits on the boundaries of Liouville domains. The construction is based on the symplectic area functional, constrained to loops with vanishing Hamiltonian mean value functional. This \emph{constrained Floer homology} (CFH) has essentially the same chain groups as Rabinowitz Floer homology (RFH), but avoids the use of a Lagrange multiplier $-$ which allows for a more intrinsic product structure.\\

We first show that the necessary Fredholm theory reduces to the known case of RFH. In particular, the standard Morse–Bott condition is sufficient.\\  

Then, we prove the $L^\infty$-bounds required for compactness of the moduli spaces. Controlling the non-local term that arises from differentiating along the constraint seems to require an additional assumption: The Liouville vector field must be of gradient type, a stronger Weinstein condition.
    \end{abstract}
    \end{center}
    \tableofcontents 
\thispagestyle{empty}


\section*{Introduction}
The original setting of Andreas Floer's infinite-dimensional Morse theory from the 1980s was that of a closed symplectic manifold $(W,\omega)$. In the case of $\partial W\neq\emptyset$, e.g. for Liouville domains, it is natural to consider similar theories for periodic orbits contained in the hypersurface $\Sigma \vcentcolon =\partial W$. 

From 1994 to 1999, different versions of what is now known as \emph{symplectic homology} were introduced by (among others\footnote{See \cite{SymHom} for an extensive survey of the different versions.}) Cieliebak, Floer\footnote{Posthumous.}, Hofer and Viterbo in \cite{FloerHofer}, \cite{CFH}, and \cite{Viterbo}. In this theory, homology groups arise as universal objects under direct and inverse limits over a partially ordered family of Hamiltonians.

Alternatively, after attaching a Lagrange multiplier to the mean value term $$\vspace{-0.35\baselineskip}  h :\loops{W} \nach \R,\ \gamma\auf \int\limits_0^1 H(\gamma(t))\d t $$ of Floer's original action functional, its critical points are precisely the closed characteristics contained in $\Sigma$. This gives rise to \emph{Rabinowitz Floer homology} (RFH), which was first presented by Cieliebak and Frauenfelder in \cite{ExactContact}. Together with Oancea, they showed in \cite{CFO} that both homologies are isomorphic.\vspace{0.25\baselineskip}

Thus, as established in \cite{VshapedProduct}, the product structure on (V-shaped) symplectic homology induces a product structure on RFH. Constructing an intrinsic product structure on RFH however is difficult, since its essential feature $-$ the Lagrange multiplier $-$ is not additive under the concatenation of loops.\vspace{0.25\baselineskip}

To circumvent this issue, we will explore yet another approach:  
Instead of using a Lagrange  multiplier, we restrict the symplectic area functional $\act{}$ to $\H\vcentcolon=h\inv(0)\subset\loops{W}$. This way, it stays additive under concatenation and still shares its critical points with the Rabinowitz action functional (\aref{thm4.8}).\vspace{0.25\baselineskip}

The main difficulty with this approach is that the gradient flow equation differs qualitatively from that of RFH: Due to the restriction to $\H$, the non-local factor \vspace{-0.2 \baselineskip} $$\mult(\gamma)= -\frac{\d\act{}_\gamma (\nabla H|_\gamma)}{\d h_\gamma(\nabla H|_\gamma)}\vspace{-0.2 \baselineskip}$$ appears. At critical points, it coincides with the period functional (\aref{thmHam}), but without further assumptions on the Liouville domain it is hard to bound:

\begin{thm} If $H\in C^\infty(M,\R)$ satisfies $\iota_{\nabla H}\d\lambda = -\lambda$,  then \vspace{-0.2 \baselineskip}
\begin{equation*}
\max\limits_{s\in\R} |\mult(u(s))| \leq C(x^\pm,H)\vspace{-0.2 \baselineskip}
\end{equation*} holds for all $u\in\M(x^-,x^+)$ and a constant $C$ only dependent on $x^-,x^+$ and $H$.
\end{thm}
This result, which is proved as \aref{lemBoundConstraint}, can be seen as an analogue to Proposition 3.2 in \cite{ExactContact}: There, similar bounds are established for the Lagrange multiplier, but without assuming $\nabla H$ to be a Liouville vector field.  On the other hand, the Fredholm theory works without any additional challenges: Using the tangent space decomposition \labelcref{eqSplitting2}, we show the following in \aref{thmFred}: 

\begin{thm} The linearization of the gradient flow operator $u\auf\partial_s u +\nabla\act{H}(u)$ is a Fredholm operator, if and only if the corresponding operator from RFH is. In this case, their indices coincide.
\end{thm}\vspace{-0.2\baselineskip}

In  \cref{secHom}, we thus may assume that the moduli spaces $\M(x^-,x^+)$ are finite-dimensional  and compact up to breaking. In other words,  the \emph{mean value constrained Floer homology} (CFH) can be constructed along the lines of the usual Morse-Bott approach known from RFH. In the subsequent paper \cite{Prod}, we will explore the product structure on CFH.\vspace{0.25\baselineskip}

\textit{Acknowledgment.} I would like to thank Prof. Urs Frauenfelder for suggesting the topic and for many helpful discussions. I am also grateful to Frederic Wagner and Zhen Gao for proofreading, and to Yannis Bähni for sharing early access to his forthcoming book \emph{Lectures on Twisted Rabinowitz–Floer Homology}.
\thispagestyle{empty}

\newpage \pagenumbering{arabic}

\section{Fredholm Theory for Moduli Spaces}\label{Chapter1}
In this section, we establish that our moduli spaces of interest are given as zeros of a Fredholm map. Under the assumption of transversality, this means that they are smooth manifolds of finite dimension.
\subsection{The Infinite Dimensional Setup}\label{SecInfDimSetup}
Let $(M,\omega)$ be an exact symplectic manifold together with a time-independent Hamiltonian $H\in C^\infty(M,\R)$. Writing $\S\vcentcolon =\R/\Z$, the \emph{loop space}  of $M$ is $\loops{M} \vcentcolon = C^\infty(\S,M)$.  
 
\definition\label{defFuncs} Let $\omega =\d\lambda$. The \emph{(negative) symplectic area functional} is given by
\begin{align*}
\act{} : \loops{M} \nach \R,\text{ } &\gamma\auf -\int\limits_{\S} \gamma\star\lambda,
\intertext{and the \emph{(Hamiltonian) mean value functional} by} 
h : \loops{M} \nach \R, \text{ }&\gamma\auf \int\limits_0^1 H(\gamma(t))\d t. 
\end{align*} 

To fix notation, we state the following standard definition and result:
\definition An \emph{almost complex structure} on $(M,\omega)$ is a smooth section\linebreak $J\in\Gamma\ndo{\tang{}{M}}$ with the property $J^2 = -\id$. It is called \emph{$\omega$-compatible}, if $\omega(\cdot,J\cdot)$ is a Riemannian metric on $M$.

\begin{lemma}\label{lemGrad} There exists a unique vector field $X_H\in\vec{M}$ with $\d H = \omega(\cdot,X_H)$. Similarly, for any Riemannian metric $g$ on $M$, there exists a unique vector field $\nabla H\in\vec{M}$ with $\d H = g(\nabla H,\cdot)$. Further, if $g=\omega(\cdot,J\cdot)$ for some almost complex structure $J$, then $\nabla H = -JX_H$.\end{lemma}

Since $\loops{M}$ is only a Fréchet space, we trade smoothness for completeness:

\theorem\label{thmBanach} For $p> 1$, the \emph{weak loop space}  $\loops[p]{M} \vcentcolon = W^{1,p}(\S,M)$ is a Banach manifold. Its tangent space at $\gamma\in\loops[p]{M}$ is given by the Banach space
\begin{equation*}
\tang{\gamma}{\loops[p]{M}} = W^{1,p}(\S,\gamma\star\tang{}{M})  =\menge{X \in W^{1,p}(\S,\tang{}{M}) : X(t)\in\tang{\gamma(t)}{M} \text{ for all }t\in\S}.
\end{equation*} 
\proof  This is a standard result in Floer and Morse theory: 
It relies on the fact that at every point $x\in M$ of a Riemannian manifold $(M,g)$, there exists a $r(x)>0$ called the \emph{injectivity radius}, as well as an open neighborhood $U\subset M$ of $x$, so that $\exp_x : B_{r(x)}(0)\subset\tang{x}{M}\nach U$ is a bijection. So for $\gamma\in\loops[p]{M}$, any $X$ in \begin{equation}\label{eqDefCharts}V_\gamma\vcentcolon =\{X\in W^{1,p}(\S,\gamma\star\tang{}{M}) : X(t)\in B_{r(\gamma(t))}(0) \text{ for all } t\in\S\}\end{equation} gives rise to a unique loop via $ X\auf \exp_\gamma X$ (in particular,  $0\auf\gamma$). This map is well defined for $p>1$: By Morrey's theorem\footnote{Cf. \cite{Brezis}, Theorem 9.12.}, each equivalence class in $W^{1,p}(\S,M)$ has a representative in $C^0(\S,M)$. The local charts\footnote{For more details like compatibility, see e.g. \cite{SchwarzMorse}, Appendix A.} of $\loops[p]{M}$ are then given by $(V_\gamma,\exp_\gamma)$, and from this the characterization of the tangent space also follows.\qed

Despite being defined for smooth loops, $\act{}$ and $h$ still work in the weak setting:

\lemma\label{lemDiff}  For $p>1$, the two functionals from \aref{defFuncs} are well-defined on $\loops[p]{M}$, 
and their differentials at $\gamma\in\loops[p]{M}$ are given by
\begin{align*}
\d \act{}_\gamma : \tang{\gamma}{\loops[p]{M}}\nach\R,\text{ } &X\auf \int\limits_0^1 \omega_{\gamma(t)}(\partial_t \gamma(t),X(t))\d t,\\
\d h_\gamma : \tang{\gamma}{\loops[p]{M}}\nach\R,\text{ } &X\auf \int\limits_0^1 \d H_{\gamma(t)}(X(t))\d t. 
\end{align*} 
\proof  See \aref{secAuxProofs}.\qed

\notation When the specific choice of $p>1$ from   \aref{thmBanach} and  \aref{lemDiff} is not important, we will also use $\loops{M}$ to denote the weak loop space $\loops[p]{M}$.

\lemma\label{lemReg} If $0$ is a \emph{regular value} of $H$, i.e. if $\d H_x$ is surjective at all $x\in H\inv(0)$, then $\d h_\gamma(\nabla H|_\gamma)>0$ for all $\gamma\in h\inv(0)$. In particular, $0$ is also a regular value of $h$.
\proof If $h(\gamma)=0$, then by the intermediate value theorem\footnote{Like in the proof of   \aref{thmBanach}, $\gamma\in\loops[p]{M}$ is assumed to be continuous.} for $t\auf H(\gamma(t))$, there has to exist at least one $t_0\in\S$ with $H(\gamma(t_0))=0$. Surjectivity of $\d H_{\gamma(t_0)}$ translates into  $\nabla H(\gamma(t_0))\neq 0$ due to   \aref{lemGrad}, and so continuity yields an $\eps>0$ as well as a constant $C>0$, so that $\norm{\nabla H(\gamma(t))}^2\geq C> 0$ for all $t\in (t_0-\eps,t_0+\eps)$, implying:
\begin{align*}
\d h_\gamma(\nabla H|_\gamma) &= \int\limits_0^1 \d H_{\gamma(t)} (\nabla H(\gamma(t)))\d t\\
&= \int\limits_0^1 \underbrace{g(\nabla H(\gamma(t)),\nabla H(\gamma(t)))}_{\geq 0}\d t\\
&\geq \int\limits_{t_0-\eps}^{t_0+\eps} \underbrace{\norm{\nabla H(\gamma(t))}^2}_{\geq C}\d t\\
&\geq 2C\eps \\
&>0
\end{align*}\nopagebreak
This already shows surjectivity, since $\d h_\gamma$ is linear and $\dim(\im\d h_\gamma) = 1$.\qed

\remark For the rest of the section, we assume that $0$ is a regular value of $H$.

\theorem\label{thmSub} The \emph{constrained loop space} $\H \vcentcolon = h\inv(0)$ is a Banach submanifold of $\loops{M}$ with codimension $1$. Its tangent space at $\gamma\in\H$ is given by the Banach space
\begin{equation*}
    \tang{\gamma}{\H} = \ker\d h_\gamma = \{X\in\tang{\gamma}{\loops{M}} : \d h_\gamma(X) = 0\}.
\end{equation*}
\proof \aref{lemReg} implies $\d h_\gamma(\nabla H|_\gamma)\neq 0$ for all $\gamma\in\H$, so \aref{thmBanachSub} applies.\qed

\corollary\label{cor4.3} At every $\gamma\in\H$, the following decomposition holds:
\begin{align}\label{eqSplitting2}\nonumber T_\gamma \loops{M} &= T_\gamma\H\oplus\langle\nabla H|_\gamma\rangle\\
&\cong T_\gamma\H\times\R\end{align}
\proof  Follows immediately from  \aref{lemReg} and (the proof of) \aref{thmBanachSub}.\qed

\definition \begin{enumerate}
    \item  The \emph{(mean value)-constrained action functional} is defined as $$\act{H} \vcentcolon = \act{}|_\H.$$
    \item At any $\gamma\in\loops{M}$, the \emph{$L^2$-scalar product} on $\tang{\gamma}{\loops{M}}$ is given by
   $$\langle \cdot,\cdot\rangle  \vcentcolon = \int\limits_{\S} g_{\gamma(t)}(\cdot,\cdot)\d t.$$
If $\gamma\in\H$, then the \emph{induced scalar product} on $\tang{\gamma}{\H}$ is given by
    $$\langle\cdot,\cdot\rangle_\H \vcentcolon = \langle \iota \cdot,\iota\cdot\rangle,$$
where $\iota:\tang{\gamma}{\H}\nach\tang{\gamma}{\loops{M}}$ is the canonical identification $w\auf w$.
    \item The \emph{gradient} $\nabla\act{H}(\gamma)$ of $\act{H}$ at $\gamma\in\H$ is defined by $$\d\act{H}_\gamma=\langle\nabla\act{H}(\gamma),\cdot\rangle_\H.$$
\end{enumerate}
\remark The gradient of $\act{}$ at some $\gamma\in\H$ will usually lie in $\tang{\gamma}{\loops{M}}$, so we have to subtract \glqq the part orthogonal to $\H$\grqq{} to get $\nabla\act{H}(\gamma)\in\tang{\gamma}{\H}$: 
\lemma Let $\gamma\in\H$, then \begin{equation}\label{grada} \nabla\act{H}(\gamma) = J(\gamma)\partial_t \gamma +\mult(\gamma) \nabla H|_\gamma\end{equation} holds with the \emph{constraint factor}\begin{equation}\label{defMult}\mult(\gamma) \vcentcolon= -\dfrac{\d\act{}_\gamma( \nabla H|_\gamma)}{\d h_\gamma( \nabla H|_\gamma)}.\end{equation}
\proof From \aref{lemDiff} and $g=\omega(\cdot,J\cdot)$, it is immediate that $\nabla\act{}(\gamma) = J(\gamma)\partial_t \gamma$ and $\nabla h(\gamma)= \nabla H|_\gamma$, so the result follows from \aref{thmGradBanach}.\qed

Now we can state the main motivation for considering $\act{H}$ in the first place: It allows us to find (smooth) periodic orbits of $X_H$ that live inside the zero-set of $H$.

\definition\label{defCrit} The \emph{critical points} of $\act{H}$ are defined as \begin{align*}
    \crit{\act{H}} &\vcentcolon = \menge{\gamma\in\H : \d\act{H}_\gamma = 0}\\
    &\hphantom{:}= \menge{\gamma\in\H : \nabla\act{H}(\gamma)=0}.
\end{align*} 

\theorem\label{thmHam} A loop $\gamma\in\H$ is a critical point of $\act{H}$, if and only if $\gamma \in C^\infty(\S,M)$  is a solution of \begin{equation}\label{eqHam}\partial_t\gamma(t)=\mult(\gamma) X_H(\gamma(t))\end{equation} with $\gamma(t)\in H\inv(0)$ for all $t\in\S$.

\proof  With $\nabla H = -JX_H$ from   \aref{lemGrad}, it is evident that in the case of $\nabla\act{H}(\gamma)=0$, \labelcref{grada} is equivalent to \labelcref{eqHam}. However,  $\gamma$ a priori only has $W^{1,p}$-regularity.

To prove smoothness, we use a technique known as \emph{bootstrapping}: As was already used in the proofs of   \aref{thmBanach} and   \aref{lemDiff}, $\gamma$ can be assumed to be continuous. This in turn means that $$t\auf \mult(\gamma)X_H(\gamma(t))$$ is also continuous, so because of \labelcref{eqHam} this property also follows for $t\auf \partial_t\gamma(t)$.\newline
In other words, we have shown that $\gamma\in W^{2,p}(\S,M)$. Repeating this reasoning inductively, we see that $\gamma$ has arbitrarily high  (weak) regularity and so  the Sobolev embedding theorem guarantees even classical smoothness.

Now, \labelcref{eqHam} implies that $t\auf H(\gamma(t))$ is constant due to
\begin{align*}
    \ddt{t} H(\gamma(t)) &= \d H_{\gamma(t)}(\partial_t\gamma(t))\\
    &= \omega_{\gamma(t)}(\partial_t\gamma(t),X_H(\gamma(t)))\\
    &\overset{\labelcref{eqHam}}{=} \mult(\gamma)\omega_{\gamma(t)}(X_H(\gamma(t)),X_H(\gamma(t)))\\
    &= 0,
\end{align*}
and because of $$h(\gamma)=\int\limits_0^1 H(\gamma(t))\d t=0,$$ this constant has to be $0$.\qed


\subsection{The Rabinowitz Action Functional}\label{secRab}
In \cite{ExactContact}, Cieliebak and Frauenfelder introduced the \emph{Rabinowitz action functional}
\begin{equation*}
    \Act{H} : \loops{M}\times \R\nach\R, \text{ } (\gamma,\lag) \auf \act{} (\gamma) + \lag h(\gamma).
\end{equation*}
Its differential is given by
\begin{equation}\label{eq16}
    \d\Act{H}_{(\gamma,\lag)} : T_\gamma\loops{M} \times\R\nach\R, \text{ } (X,\l) \auf \d\act{}_\gamma(X) + \lag\d h_\gamma(X)  + h(\gamma)\l, 
\end{equation} 
\begin{align}
    \intertext{and its gradient with respect to the scalar product }\nonumber \langle(X_1,\l_1),(X_2,\l_2)\rangle_{(\gamma,\lag)} &\vcentcolon = \langle X_1,X_2\rangle_{\gamma}+\l_1\l_2
    \intertext{on $\tang{\gamma}{\loops{M}}\times\R$ is given by}
   \label{gradA} \nabla\Act{H}(\gamma,\lag) &\hphantom{:}= \left(\begin{matrix}J(\gamma)\partial_t \gamma +\lag \nabla H|_\gamma\\ h(\gamma) \end{matrix}\right).
\end{align}

It is connected to the constrained action functional from the previous subsection in several important aspects, the first being the following:

\theorem\label{thm4.8} $(\gamma,\lag)\in \crit {\Act{H}}$, if and only if $\gamma\in\crit{\act{H}}$ and $\lag = \mult(\gamma)$.
\proof 
The second component of \labelcref{gradA} clearly vanishes if and only if $\gamma\in\H$, so we can assume $h(\gamma)=0$ and only have to focus on the first component:

\glqq$\Leftarrow$\grqq{}: If $\nabla\act{H}(\gamma)=0$, then $J(\gamma)\partial_t\gamma +\lag\nabla H|_\gamma$ obviously vanishes for $\lag = \mult(\gamma)$.

\glqq$\Rightarrow$\grqq{}: If $\d\Act{H}_{(\gamma,\lag)}=0$, then \labelcref{eq16} with $X=\nabla H|_\gamma$ implies $\lag = -\frac{\d\act{H}_\gamma(\nabla H|_\gamma)}{\d h_\gamma(\nabla H|_\gamma)} $, and so the first component of \labelcref{gradA} implies $\nabla\act{H}(\gamma)=0$. \qed  

In particular,   \aref{thmHam} carries over:
\begin{corollary}  A pair $(\gamma,\lag)\in\loops{M}\times\R$ is a critical point of $\Act{H}$, if and only if $\gamma \in C^\infty(\S,M)$  is a solution of \begin{equation}\partial_t\gamma(t)=\lag  X_H(\gamma(t))\end{equation} with $\gamma(t)\in H\inv(0)$ for all $t\in\S$.\end{corollary}

\begin{remark} Both $\act{H}$ and $\Act{H}$ are designed to find periodic Hamiltonian orbits on $H\inv(0)$. To achieve this, they rely on different techniques: $\act{H}$ is \emph{constrained} to the hypersurface $\H\subset\loops{M}$ from the beginning and $\mult$ appears only after differentiating along this constraint, while $\Act{H}$ uses the \emph{Lagrange multiplier} $\lag$ to turn the desired constraint into part of the gradient.\end{remark}


\subsection{Moduli Spaces of Gradient Flow Cylinders}\label{secGradFlowCyl}
Next, we adapt the concept of gradient flow cylinders to the constrained setting.
\begin{definition}\label{defModuliSpaces} A  \emph{negative gradient flow cylinder} (or \emph{Floer cylinder}) is a smooth solution $u\in C^\infty(\R,\H)$ of \begin{equation}\label{eqGradFlow01} \partial_su = -\nabla \act{H}(u).\end{equation}
In other words, $u\in C^\infty(\rs,M)$ satisfies 
\begin{equation*}
    \left\{\begin{matrix}\partial_s u + J(u)\partial_t u +\mult(u)\nabla H|_u &=0\hphantom{.}\\
   h(u) 
   &=0.
   \end{matrix}\right. 
\end{equation*}
Given two critical points $x^\pm\in\crit{\act{H}}$, the \emph{moduli space of connecting gradient flow cylinders} is given by \begin{equation*} \M(x^-,x^+) \vcentcolon = \menge{u\in C^\infty(\R,\H) : u\text{ solves \labelcref{eqGradFlow01} and } u(s,\cdot) \overset{C^1}{\nach} x^\pm\text{ for }s\nach\pm\infty}.\end{equation*}
\end{definition} 

\begin{claim}Under certain assumptions (specified later), these moduli spaces  are finite-dimensional manifolds.\end{claim}

We will prove this \glqq standard result\grqq{} in the usual way, i.e. by constructing charts from the zero sections of Fredholm operators. As in the case of loop spaces, we first have to leave the smooth setting to get completeness:  After fixing $x^\pm\in\crit{\act{H}}$, we  consider the vector bundle $L^p\nach\W^p$ over
$$\W^p  \vcentcolon = \menge{u\in W^{1,p}(\rs,M) : u(s)\in\H\text{ for all }s\in\R\text{ and } u(s,\cdot)\underset{s\rightarrow\pm\infty}{ \overset{W^{1,p}}{\nach}}x^\pm}$$
with fiber $$L^p(u)\vcentcolon= L^p(\R,u\star\tang{}{\H})$$ 
over every $u\in\W^p$. 

\theorem\label{thmCompl}  For $p>2$, $L^p\nach \W^p$ is a Banach bundle.

\proof Like in the proof of \aref{thmBanach}, the exponential map is used to construct charts of $\W^p$, while $L^p$ naturally carries the structure of a Banach space. The main difference is that because of $\dim(\rs)=2$, the condition $p>2$ is needed to find continuous representatives. For more details, see e.g. \cite{FloerMorse}, Chapter 8.2.\qed

\definition The \emph{gradient flow operator} of $\act{H}$ is the section
\begin{align}\label{eqGradFlow03T}
\F : \W^p\nach L^p,\ u\auf \partial_s u + J(u)\partial_t u + \mult(u)\nabla H|_u,
\end{align} i.e. $\F(u)\in L^p(u)$ for every $u\in\W^p$.

\newpage
\theorem\label{thmReg} For $p>2$, any  $u\in\W^p$ with $\F(u)=0$  lies in $C^\infty(\R,\H)$. In particular, it holds that $\M(x^-,x^+)=\F\inv(0)$.
\proof Similar to   \aref{thmHam}, this time however, there are two partial derivatives and asymptotics to think about. This makes the concrete arguments much more involved, which is why we refer to \cite{FloerMorse}, Chapter 12.\qed
 
With this characterization of the moduli spaces in our minds, we now state the tools that we want to apply in order to prove the claim:
\definition A map $f\in C^1(\X ,\Y)$ between two Banach manifolds $\X$ and $\Y$ is called a \emph{Fredholm map} if its linearization $$\d f_x : \tang{x}{\X}\nach\tang{f(x)}{\Y}$$ is a Fredholm operator at every $x\in\X$, i.e. if $\d f_x$ has closed image and both its kernel and cokernel are finite dimensional. Its \emph{index} is then defined by $$\ind(f) = \dim(\ker\d f_x) - \dim(\coker\d f_x)$$ and independent of $x\in\X$.

\theorem\label{thmRegFred} Let $f : \X \nach \Y$ be as above. If $y\in \Y$ is a regular value of $f$, then $f\inv(y)$ is a submanifold of $\X$ with dimension equal to $\ind(f) $.

\proof That $f\inv(y)$ is a submanifold follows from the infinite-dimensional version of the regular value theorem also used in   \aref{thmSub}. Because of $\tang{x}{f\inv(y)}=\ker\d f_x$ and $\dim(\ker\d f_x)<\infty$ due to the Fredholm property, the dimension of $f\inv(y)$ follows from the definition of $\ind(f)$ above. \qed
\vspace{-0.5\baselineskip}
\lemma\label{lemLin} For $p>2$, the linearization of $\F$ at $u\in\W^p$ 
is the map 
$$\d\F_u :  W^{1,p}_\delta(\rs,u\star\tang{}{M})\cap L^p(\R,u\star\tang{}{\H})\nach L^p_\delta(\R,u\star\tang{}{\H})  
$$ between $\delta$-weighted Sobolev spaces\footnote{One says $f\in W^{k,p}_\delta(\rs,\R^{2n})$, if  $\gamma_\delta f\in W^{k,p}(\rs,\R^{2n})$. Here, $\gamma_\delta(s) \vcentcolon= e^{\delta\beta(s) s}$ is defined using some $\beta\in C^\infty(\R,[-1,1])$ which is strictly increasing in $[-s_0,s_0]$ and constantly $\pm 1$ outside. The importance of the weights for our setting will become apparent only in the next paper.} 
given by
\begin{align}\label{eqLin}
\xi \auf \nabla_s \xi +\nabla_\xi(J(u)\partial_t u+\mult(u)\nabla H|_u). 
\end{align}\vspace{-1\baselineskip}
\proof In   \aref{thmCompl}, the Banach space which $\W^p$ is modelled on was left ambiguous: Usually, one takes $W^{1,p}(\rs,u\star\tang{}{M})$ and then constructs all cylinders close to $u$ using the exponential map and vector fields in $W^{1,p}$ with lengths smaller than the injectivity radius at every point of $u$. But since $W^{1,p}_\delta$ differs from $W^{1,p}$ only by a global multiplication by $\gamma_\delta$, the construction goes through just the same. 
The calculation of $\d\F_u$ 
is long but standard: For the general idea, see e.g. \cite{Topol}, Section 3.3.3. A proof in coordinates is given in \cite{FloerMorse}, Section 8.4.\qed

\newpage

\subsection{Completions of Liouville Domains}\label{secSympSetup}
Before proceeding with the Fredholm theory, let us fix the type of manifold we are interested in and define a canonical form of $H$.
\definition Let $(W,\omega)$ be a compact symplectic manifold with boundary. A \emph{Liouville vector field} is a $\liou\in\vec{M}$ so that $\lambda\vcentcolon =\iota_\liou\omega$ satisfies\footnote{Clearly, the existence of a Liouville vector field is equivalent to exactness of $M$.} $\d\lambda = \omega$.  If $\liou$ is transverse to the boundary $\partial W$ and outward pointing, the tuple $(W,\liou,\lambda)$ is called a \emph{Liouville domain}.

\lemma\label{lemCoords} Let $(W,\liou,\lambda)$ be a Liouville domain, and write $\Sigma \vcentcolon = \partial W$ as well as $M^+\vcentcolon =[0,\infty)\times\Sigma$ with coordinates $(r,x)$. Further, let $$M\vcentcolon =W\sqcup_\Sigma M^+$$ be the manifold obtained by identifying $\partial W$ with $\partial M^+$, as well as $\liou|_\Sigma$ with $\partial_r|_{\{0\}\times\Sigma}$. The triple $(M,\Sigma,\lambda)$ is called \emph{completion} of $(W,\liou,\lambda)$, 
Then there exists an $\eps>0$, an open set $U_\eps\subset M$ containing $M^+$, and coordinates $(-\eps,\infty)\times\Sigma \nach U_\eps$, so that 
\begin{align*}
 (-\eps,0)\times\Sigma&\inach U_\eps\cap(W\setminus\Sigma),\\
\{0\}\times\Sigma&\inach \Sigma,\\
(0,\infty)\times\Sigma&\inach M^+\setminus\Sigma.
\end{align*}Further, there exists an exact symplectic form $\bar\omega$ defined on all of $M$ that coincides with $\omega$ on $W$ and can be written in coordinates  as 
\begin{equation}\label{eqOmega}\bar\omega_{r} = e^r(\d\lambda|_\Sigma+\d r\wedge \lambda|_\Sigma)\end{equation}for all $r\in(-\eps,\infty)$. Its primitive is given by $\bar\lambda_r = e^r\lambda|_\Sigma$.

\proof Since $\liou$ is transversal to the hypersurface $\Sigma$, there exists an $\eps>0$ and an open neighborhood $U_\eps\subset M$ of $\Sigma$ so that the flow $\Phi_\liou: (-\eps,\eps)\times\Sigma\nach U_\eps$ is a diffeomorphism\footnote{This is a standard result called \emph{Flowout theorem}, see e.g. \cite{LeeSmooth}, Theorem 9.20(d).}. For  $r\geq 0$, we also have $\liou=\partial_r$ and thus $\Phi_\liou^r(x)=(r,x)$ for all $x\in\Sigma$, i.e. $\Phi_\liou$ can be extended to a diffeomorphism from $[0,\infty)\times\Sigma$ to itself. Transversality and outward-pointing of $\liou$ further implies that $(-\eps,0)\times\Sigma$ is contained in $W\setminus\Sigma $. Next, for $(r,x)\in(-\eps,0]\times\Sigma$,  observe that 
\begin{align}
\nonumber\tilde\lambda_{(r,x)} &\vcentcolon=(\Phi^r_\liou(x))\star\lambda\\
\nonumber&= e^r\lambda_{\Phi^0_\liou(x)}\\
&= e^r\lambda_{x},\label{eqLambdaTilde}
 \end{align}
where the second equality follows from a short calculation involving the Lie derivative\footnote{See e.g. the proof of (2.17) in \cite{ThreeBody}.}. $\tilde\lambda$  can be extended to all $r\in(-\eps,\infty)$ using \labelcref{eqLambdaTilde}, so   $\bar\omega\vcentcolon=(\Phi^{-r}_\liou)\star\d\tilde\lambda$ is defined on all of $U_\eps$. Finally, pullback and differential commute, so $\bar\omega$ agrees with $\omega$ on $U_\eps\cap W$ due to the definition of $\tilde\lambda$ and can be extended to all of $W$ using $\omega$.\qed

\remark In the following, we will neither distinguish between the symplectic form $\omega$ on $W$ and its extension $\bar\omega$ to all of $M$, nor between $\lambda$ and $\bar\lambda$.

\definition The \emph{Reeb flow} of a contact manifold $(\Sigma,\alpha)$ is implicitly defined by $\iota_R\alpha = 1$ and $\iota_R\d\alpha = 0$. If $\Sigma\subset M$ is a hypersurface inside a symplectic manifold $(M,\omega)$ with $\d\alpha=\omega|_\Sigma$, then $H\in C^\infty(M,\R)$ is called a \emph{defining Hamiltonian} for $(\Sigma,\alpha)$ if it has $0$ as a regular value, $\Sigma = H\inv(0)$, and $X_H=R$.

\lemma\label{lem1.29} If $(M,\Sigma,\lambda)$ is the completion of a Liouville domain, then $(\Sigma,\lambda|_\Sigma)$ is a contact manifold. 
A defining Hamiltonian for $\Sigma$ is given by any $H\in C^\infty(M,\R)$ with $\Sigma = H\inv(0)$ and $$H(r,x)= e^r-1$$ in an open neighborhood around $\Sigma$ (and using the coordinates from \aref{lemCoords}).
\proof
 The first property is a standard result, see e.g. Lemma 2.6.3 in \cite{ThreeBody}. For the second, note that because of $\d H = e^r \d r$, the condition $\d H = \omega(\cdot,X_H)$ with $\omega$ as in \labelcref{eqOmega} immediately implies that $\iota_{X_H}\d\lambda |_\Sigma =0$ and $\iota_{X_H}\lambda|_\Sigma =1$.
\qed 

\definition\label{defSFTlikeJ} Let $(M,\Sigma,\lambda)$ be the completion of a Liouville domain and let $\xi\vcentcolon =\ker\lambda|_\Sigma$ be the contact distribution of the contact manifold $(\Sigma,\lambda|_\Sigma)$. A $\d\lambda$-compatible almost complex structure $J$ on $M$ is called \emph{SFT-like} if $J\liou=R$, 
$J|_\xi$ is $\d\lambda|_\xi$-compatible, and $J_{(r,x)}=J_{(0,x)}$ for all $(r,x)\in (-\eps,\infty)\times\Sigma$. 

\begin{remark} Clearly, any $\d\lambda|_\xi$-compatible almost complex structure $J_\xi$ on $\xi$ can be extended to an SFT-like almost complex structure $J$ on $M$. 
\end{remark}

\corollary\label{corGrad}  The gradient of any defining Hamiltonian $H$ as in \aref{lem1.29} with respect to the scalar product  induced by an SFT-like almost complex structure satisfies $$\nabla H = \liou$$ in a neighborhood around $\Sigma$. In particular, \begin{equation}\label{eqNormNablaH}\norm{\nabla H(r,\cdot)}^2 = \norm{\liou(r,\cdot)}^2=e^r\end{equation} holds for all  $r\geq -r_0$ with $r_0\in (0,\eps)$.
\proof Obvious from $\d H = e^r\d r$ and \aref{defSFTlikeJ}.\qed


\subsection{The Fredholm Property}
We will now state and prove the main result of this section: First off, all the concepts from \cref{secGradFlowCyl} can be transferred to the Rabinowitz action functional $\Act{H}$: Its gradient flow cylinders $(u,\lag)\in C^\infty(\R,\loops{M})\times C^\infty(\R,\R)$ are the solutions to
\begin{equation*} \partial_s (u,\lag) = -\nabla \Act{H}(u,\lag),\end{equation*} i.e. the $(u,\lag)\in C^\infty(\rs,M)\times C^\infty(\R,\R)$ with
\begin{equation*} 
    \left\{\begin{matrix}\partial_s u + J(u)\partial_t u +\lag \nabla H|_u &=0\hphantom{.}\\
   \partial_s\lag + h(u) 
   &=0,
   \end{matrix}\right. 
\end{equation*}
and its gradient flow operator is a section $\G$ in a Banach bundle similar to \labelcref{eqGradFlow03T}. From \cite{ExactContact}, it is known that its linearization
$$\d\G_{(u,\lag)} : W^{1,p}_\delta(\rs, u\star\tang{}{M})\times W^{1,p}_\delta(\R,\R) \nach L^p_\delta(\rs,u\star \tang{}{M})\times L^p_\delta(\R,\R)$$at $(u,\lag)$ has the form
\begin{equation}\label{eq20}
    (X,\l) \auf \left(\begin{matrix}\nabla_s X + \nabla_X( J(u) \partial_t u+\lag\nabla H|_u) + \l\nabla H|_u\\
    \l' + \d h_u(X)\end{matrix}\right) 
\end{equation}
and is a Fredholm operator for a suitable $\delta>0$ under the following assumption:

\assumption\label{assSymp} Let $(M,\Sigma,\lambda)$ be the completion of a Liouville domain, and assume that $H\in C^\infty(M,\R)$ is a defining Hamiltonian of $(\Sigma,\lambda|_\Sigma)$ like in \aref{lem1.29}. Further, let $ \act{H} $ be \emph{Morse-Bott}, i.e. $\crit{\act{H}}\subset\Sigma$ is a properly embedded submanifold (see also \textup{Assumption (A)} in \textup{\cite{ExactContact}}, p. 273).

\remark For every $\gamma\in\crit{\act{H}}$ and $r\in\S$, another critical point is given by $r_*\gamma\vcentcolon=\gamma(\cdot+r)$ and so $\act{H}$ can never be an actual Morse function. In other words, its \emph{Hessian} $\hess{\act{H}}$ has $0$ as an eigenvalue. Under the Morse-Bott condition, and for $0<\delta<\min\{|\lambda| : \lambda \in\spec(\hess (\act{H}))\setminus\{0\}\}$, $\Act{H}$ nonetheless forms a Fredholm operator between the $\delta$-weighted Sobolev spaces above.

\begin{theorem}\label{thmFred}
    Under \aref{assSymp} and for a suitable $\delta>0$, the linearization $\d\F_u$ of $\F$ at $u\in\W$  
    is a Fredholm operator with the same index as $\d\G_{(u,\mult(u))}$.
\end{theorem}

\proof Since $u\in\W$, we can use the decomposition \labelcref{eqSplitting2} of $\tang{u(s)}\loops{M}$ into $\tang{u(s)}{\H}$ and $\langle\nabla H|_u(s)\rangle$. So for any $(X,\l)\in W^{1,p}(\rs,u\star\tang{}{M})\times W^{1,p}(\R,\R)$, we can write
$$X(s)=\xi(s)\oplus\alpha(s)\nabla H|_u(s).$$

With this notation, the first component of \labelcref{eq20} takes the form
\begin{equation*} \underbrace{\nabla_s\xi +\nabla_\xi(J(u)\partial_tu+\mult(u)\nabla H|_u)}_{=\d\F_u(\xi)}+ \nabla_s \alpha\nabla H|_u + \nabla_{\alpha\nabla H|_u}(J(u)\partial_t u+\mult(u)\nabla H|_u)+\l\nabla H|_u \end{equation*}
and we immediately identify the first part as the operator $\d\F_u$ from \labelcref{eqLin} which we are interested in. To relate it to $\d\G_{(u,\mult(u))}$, we will now show that all the other terms do not contribute to the Fredholm properties.

First off, note that $\lim\limits_{s\nach\pm\infty} u(s,\cdot)=x^\pm \in\Sigma,$ so for sufficiently large $s\in\R$, the cylinder has to lie in the open set $U_\eps$ from \aref{lemCoords} containing $\Sigma$. Thus, after a small homotopy, we may assume that $\im(u)\subset U_\eps$. Since $\d\G_{(u,\lag)}$ remains a Fredholm operator along this smooth family of cylinders, its Fredholm index does not change.  

First off, using elementary properties of the covariant derivative, we see
\begin{align*}
    \nabla_s\alpha\nabla H|_u &= \alpha'\nabla H|_u + \alpha \nabla_{\partial_s u}\nabla H|_u\\
     \shortintertext{as well as}
 \nabla_{\alpha\nabla H|_u} (\cdots) &=\alpha\nabla_{\nabla H|_u} (\cdots).
\end{align*}

Using these results and $\nabla_{\partial_s u}\nabla H|_u = \nabla_{\nabla H|_u}\partial_s u$, we can write the remaining terms of the first component as follows:
\begin{align*} \alpha(s)\nabla_{\nabla H|_u} (\partial_s u + J(u)\partial_t u+\mult(u)\nabla H|_u) + (\alpha'(s)+\l)\nabla H|_u
\end{align*}

For the second component, we note 
$$\d h_u(X)=\alpha \d h_u(\nabla H|_u).$$ In total, this means that we can write $\d\G_{(u,\mult(u))}$ from \labelcref{eq20} as
\begin{align*}
W^{1,p}_\delta(\rs,u\star\tang{}{M}) \cap L^p (\R,&u\star\tang{}{\H})\times W^{1,p}_\delta(\R,\R^2) \nach L^p_\delta(\R,u\star\tang{}{\H})\times L^p_\delta(\R,\R^2),\\
(\xi,\alpha,\l)&\auf \left(\begin{matrix} \d\F_u(\xi)\\ A(\alpha,\l)\end{matrix}\right)+K(\alpha,\l). 
\shortintertext{Here, we used $\langle \nabla H|_u(s)\rangle \cong\R$ and the two auxiliary operators}
K: W^{1,p}_\delta(\R,\R)&\nach W^{1,p}_\delta(\rs,u\star\tang{}{M})\hookrightarrow L^p_\delta(\rs,u\star\tang{}{M}),\\
(\alpha,\l)&\auf  \alpha \nabla_{\nabla H|_u}(\partial_s u + J(u)\partial_t u +\mult(u)\nabla H|_u)  
\intertext{and}
A : W^{1,p}_\delta(\R,\R^2)&\nach L^p_\delta(\R,\R^2),\\
(\alpha,\l)&\auf  
\left(\begin{matrix} \alpha' +\l \\ \d h_u(\nabla H|_u)\alpha+ \l'\end{matrix}\right).
\end{align*} 

If $u$ was a gradient flow cylinder (i.e. $\G(u,\mult(u))=0$), $K$ would simply vanish. In the case that it does not, we can assume after another homotopy that there exists some $T>0$ with $u(s,\cdot)=x^\pm$ for $|s|>T$. Then, $(K\alpha)(s)=0$ for all $|s|>T$, so $K$ can be seen as the composition  
$$K:W^{1,p}_\delta(\R,\R)\overset{}{\nach} W^{1,p}_\delta([-T,T],\R) \overset{\iota}{\hookrightarrow} L^p_\delta(\rs,u\star\tang{}{M})$$ which is clearly compact and thus can be neglected from the Fredholm analysis.

Further, $ (\d h_{u}(\nabla H|_u))(s) = 1$ for $|s|>T$, so $$A  = \partial_s + \left(\begin{matrix}0 & 1\\ \d h_u(\nabla H|_u) & 0 \end{matrix}\right)$$ is a Fredholm operator well-known from Morse theory: Its index is zero, since the asymptotic matrices for $s\nach\pm\infty$ coincide  (see e.g. \cite{Schwarz}, Propositon 2.16). 

Thus, we can write $\d\G_{(u,\mult(u))}$ as the direct sum of $\d\F_u$ and $A$, which immediately implies that $\d\F_u$ has to be Fredholm too and that $$\ind(\d\G_{(u,\mult(u))})=\ind(\d\F_u)+\ind(A)=\ind(\d\F_u). $$ \qed

\begin{corollary}\label{corMF}
As long as \aref{assSymp} holds and $\d\F_u$ is surjective at all $u\in\F\inv(0)$, the moduli space $\M(x^-,x^+)$ is a manifold of $\dim\M(x^-,x^+)=\ind(\F)$.
\end{corollary}
\proof Because of \aref{thmFred} and surjectivity, \aref{thmRegFred} can be applied.\qed

\begin{remark}
 It is known from Theorem B.1 in \textup{\cite{ExactContact}} that $\d\G_{(u,\tau)}$ is surjective at all $(u,\tau)\in\G\inv(0)$ for a \emph{generic} choice of defining Hamiltonian $H$. To achieve the manifold-property without perturbing $H$ or $J$, one can also pass to an \emph{abstract perturbation}  in the sense of polyfold theory: See e.g. \textup{Chapter 5} in \textup{\cite{YannisBook}}.
\end{remark}

\newpage

\section{Compactness of Moduli Spaces}\label{secCompMod}
In this section, we always suppose \aref{assSymp} holds.
\subsection{A Bound on the Constraint Factor}\label{secBoundConstraint}
To achieve a mean value of zero, i.e. $$\int\limits_0^1 H(\gamma(t))\d t=0,$$ it is clear that every $\gamma\in\H$ has to spend some time in a region of $M$ where $\nabla H>0$, as $0$ is a regular value of $H$ However, it might be possible that a sequence $(\gamma_n)_{n\in\N}\subset\H$ exists for which $\d h_{\gamma_n} (\nabla H|_{\gamma_n})$ tends to zero. In the following lemma, we will show that this in fact cannot happen by providing a uniform version of \aref{lemReg}:
\lemma\label{lemBoundDiff} For any $H$ as in \aref{lem1.29}, there exists a constant $c>0$, so that
$$\d h_\gamma(\nabla H|_\gamma)>c\ \forall\gamma\in\H.$$ 

\proof For arbitrary $\gamma\in\H$, we can define (using the coordinates $(r,x)$ from \aref{lemCoords})
\begin{align*}
G^+ &\vcentcolon=\menge{ t\in [0,1] : H(\gamma(t))\geq 0},\\
    G^- &\vcentcolon= \menge{ t\in [0,1] : H(\gamma(t))< 0},\\
    \widetilde G\hphantom{^-} &\vcentcolon= \menge{t\in G^- : r(\gamma(t))>-r_0}.
\shortintertext{From the fact that $h(\gamma)=0$, we immediately get}
\int\limits_{G^+} H(\gamma(t))\d t &= -\int\limits_{G^-} H(\gamma(t))\d t,
\shortintertext{and from \aref{corGrad}, we already know that $\norm{\nabla H(r,\cdot)}^2 =e^r$ for $r>-r_0$. Thus,}
\int\limits_{G^+}\norm{ \nabla H(\gamma(t))}^2\d t &= \hphantom{-}\int\limits_{G^+} e^{r(\gamma(t))}\d t\\
&= \hphantom{-}\int\limits_{G^+} H(\gamma(t))+ 1\d t\\
&= -\int\limits_{G^-} H(\gamma(t))\d t + \int\limits_{G^+} 1\d t\\
&= -\int\limits_{\widetilde G} \underbrace{H(\gamma(t))}_{\geq 0} \d t - \int\limits_{G^-\setminus \widetilde G} \underbrace{H(\gamma(t))}_{\leq -a}\d t + \mu( G^+)\\
&\geq a\mu(G^-\setminus \widetilde G) + \mu(G^+)
\shortintertext{using the Lebesgue measure $\mu$ and a constant}
 a &\vcentcolon= -\max \menge{H(z) : z\in W\setminus \{(r,x) \in U_\eps : r >-r_0\}}
 \shortintertext{which is independent of $\gamma$ and positive: $M$ and $\Sigma =H\inv(0)$ are both connected, so $W\setminus\Sigma =H\inv( (-\infty,0))$ and because $W\setminus \{z \in W : r(z)>-r_0\}$ is compact, $H$  attains a maximum there which is strictly smaller than $0$. Similarly, we have}
\int\limits_{G^-} \norm{\nabla H (\gamma(t))}^2\d t &= \int\limits_{\widetilde G} \norm{\nabla H (\gamma(t))}^2\d t + \int\limits_{G^-\setminus\widetilde G} \norm{\nabla H (\gamma(t))}^2\d t\\
 &\geq b\mu(\widetilde G),
 \shortintertext{for another constant}
 b &\vcentcolon= \min\menge{\norm{\nabla H(r,\cdot)} : r>-r_0}\\
 &= e^{-r_0}\\
 &>0.
 \shortintertext{Putting all of this together, we can estimate}
\d h_\gamma(\nabla H|_\gamma) &= \int\limits_0^1 \norm{ \nabla H(\gamma(t))}^2\d t\\
&= \int\limits_{G^+} \norm{ \nabla H(\gamma(t))}^2\d t + \int\limits_{G^-} \norm{ \nabla H(\gamma(t))}^2\d t \\
 &\geq a\mu(G^-\setminus \widetilde G) + \mu(G^+) + b\mu(\widetilde G)\\
 &\geq c\mu([0,1])\\
 &= c
\end{align*}for $c \vcentcolon= \min\menge{a,b,1}>0$, since all three sets are disjoint.\qed 

\lemma\label{lemBoundConstraint} If $H\in C^\infty(M,\R)$ satisfies $\nabla H =\liou $, then  
$|\mult | \leq \frac{1}{c} |\act{H} |.$ In particular,
\begin{equation}\label{eqBoundConstraint}
\max\limits_{s\in\R} |\mult(u(s))| \leq C(x^\pm,H)
\end{equation}holds for all $u\in\M(x^-,x^+)$ and a constant $C$ only dependent on $x^-,x^+$ and $H$. 

\proof The first part follows from the definition of $\mult$ in \labelcref{defMult}, the fact that $\d\act{H}_\gamma(\liou|_\gamma)=\act{H}(\gamma)$, and \aref{lemBoundDiff}. 
For \labelcref{eqBoundConstraint}, recall that $\act{H}(x^-)\geq \act{H}(u(s))\geq \act{H}(x^+)$ for $s\in\R$ and $u\in\M(x^-,x^+)$.\qed
 
\vspace{-0.75\baselineskip}
\begin{remark} The assumption $\nabla H=\liou$ is a stronger variant of the \emph{Weinstein condition}, which only requires that  $\liou$ be \emph{gradient-like} for $H$, i.e., that $\d H(\liou)>0$ away from critical points. Intuitively, it arises from the non-locality of the gradient flow of $\act{H}$\textup{:} The constraint factor $\mult$ depends on the whole loop. In contrast, the Lagrange multiplier $\lag$ in RFH only depends on $s\in\R$, and so an $L^\infty$-bound can be deduced from the gradient flow equation itself (cf. \textup{\cite{ExactContact}}, Proposition 3.2).\nopagebreak
\end{remark}


\subsection{Bounds on the Gradient Flow Cylinders}
In this subsection, we will establish $L^\infty$-bounds for (the derivatives of) gradient flow cylinders of \emph{finite energy}. This will follow from the standard result \aref{lem3.4}, once we have shown in \aref{thm3.3} that no sequence $\folge{u}{n}\subset \M(u^-,u^+)$ can be unbounded in the radial direction $[0,\infty)\times\Sigma$ of the completed Liouville domain.

\definition The \emph{energy} of a curve $u\in C^\infty(\rs,M)$ is given by
\begin{equation*}
E(u) \vcentcolon = \intlim{\R\times\S} \norm{\partial_s u(s,t)}^2\dt(s,t).
\end{equation*}
This can be simplified using the gradient flow equation:
\lemma\label{lemEnergy}  For any $u\in\M(x^-,x^+)$, its energy is finite and  can be expressed as
\begin{equation*}
E(u) = \act{H}(x^-)-\act{H}(x^+).
\end{equation*}
\proof Follows from the fundamental theorem of calculus and \labelcref{eqGradFlow01}.\qed

Our aim is to use the following standard result:
\lemma\label{lem3.4} No sequence of finite energy gradient flow cylinders with images contained in a compact subset of $M$ can have unbounded derivatives.
\proof 
This is a well-known result: A rescaling argument shows that any sequence with unbounded derivatives would yield a non-constant $J$-holomorphic sphere, which cannot exist due to the exactness of $(M,\d\lambda)$. For details, see e.g. \cite{FloerMorse}, Proposition 6.6.2, p. 176. \qed

To use \aref{lem3.4}, we only have to show that there can be no sequence  $\folge{u}{n}\subset \M(u^-,u^+)$ which is unbounded in the radial direction. For this, we first prove:

\lemma\label{lem3} Let $u\in C^\infty(\R\times\S,M)$ be any smooth map. If $u(s,t)\in [0,\infty)\times\Sigma$, we write $r(s,t)$ for the $r$-component of $u$ and have $$ \bet{\partial_t r(s,t)}\leq \norm{\partial_t u(s,t)}.$$
\proof Consider the function $r : M^+\nach [0,\infty)$, then \labelcref{eqOmega} implies $$\nabla r (z) = e^{-r(z)}\liou(z)$$ and so in general $\norm{\nabla r}^2=e^{-r}\leq 1$. Applying $r$ to the \glqq curve\grqq{} $t\auf u(s,t)$ and differentiating then yields
\begin{align*}
\partial_t r(u(s,t)) &\overset{\hphantom{\text{C.S.}}}= \d r (\partial_t u(s,t))\\
&\overset{\hphantom{\text{C.S.}}}= g_{u(s,t)}(\partial_t u(s,t),\nabla r(u(s,t)))\\
&\overset{\text{C.S.}}{\leq} \norm{\partial_t u(s,t)}\norm{\nabla r (u(s,t))}\\
&\overset{\hphantom{\text{C.S.}}}{\leq} \norm{\partial_t u(s,t)}.\qedhere
\end{align*} 

\begin{remark}
 For simplicity, we will denote the  $r$-component of $u\in C^\infty(\R\times\S,M)$ by $r \vcentcolon = r(u) \in C^\infty(\R\times\S,\R)$ and keep any indices that $u$ may have.    
\end{remark} 

The following crucial lemma is a formalization of an intuitive idea: $H$ tends to $+\infty$ for $r\nach+\infty$ and is bounded from below by some negative value, so the further some $\gamma\in\loops{M}$ with $$\int\limits_0^1 H(\gamma(t))\d t =0$$ ventures into the $r$-direction, the less time it can spend there. In other words, if the $r$-component of a sequence $\folge{\gamma}{n}\subset\H$ is unbounded, then so are their velocities.

\lemma\label{lem3.3} Let $(u_n)_{n\in\N} \subset C^\infty(\R,\H)$ be a sequence of curves, set $$R_n\vcentcolon =\max\{r_n(s,t):(s,t)\in\R\times\S\},$$  and let $(s_n,t_n)\in\rs$ be so that $r_n(s_n,t_n)=R_n$. If $\lim\limits_{n\rightarrow\infty} R_n=\infty$, then for any $\eps>0$, there exists a sequence $\folge{\theta}{n}\subset\S$  
with
\begin{align*}
  |t_n-\theta_n|&\leq \frac{\eps}{\alpha_n},\\
    |\partial_t u_n(s_n,\theta_n)| &\geq \frac{\alpha_n}{\eps},
\end{align*}
for  $\alpha_n \vcentcolon = \frac{R_n-1}{H^-}$ and all sufficiently large $n\in\N$. Here, $H^-\vcentcolon= -\min H>0 $. 
 
\proof In addition to the notation above, we define (as in the proof of \aref{lemBoundConstraint})
\begin{align*}
G^+_n &\vcentcolon =\{t\in\S : u_n(s,t) \in M^+\},\\
G^-_n &\vcentcolon =\{t\in\S : u_n(s,t)\notin M^+\},\\
\widetilde G_n &\vcentcolon = \left\{t\in G_n^+ :  r_n(s_n,t)\geq R_n-1 \right\}.
\end{align*}
As a direct consequence of  $h(u_n(s_n,\cdot))=0$, one again has
\begin{equation}\label{eq6}
\intlim{G_n^+} H(u_n(s_n,t))\d t = -\intlim{G_n^-} H(u_n(s_n,t))\d t.
\end{equation}
Also, for $t\in \widetilde G_n$, it holds that 
\begin{align*}
e^{R_n-1}-1&\leq  e^{r_n(s_n,t)}-1\\
&= H(u_n(s_n,t)),
\intertext{so for any $\eps>0$ and all sufficiently large $n\in\N$, we get}
\frac{1}{\eps}(R_n-1)\mu(\widetilde G_n) &\overset{\hphantom{\labelcref{eq6}}}{\leq} \rund{e^{R_n-1}-1}\mu(\widetilde G_n)\\
 &\overset{\hphantom{\labelcref{eq6}}}{\leq} \intlim{\widetilde G_n} H(u_n(t))\d t\\
  &\overset{\hphantom{\labelcref{eq6}}}{\leq} \intlim{G^+_n}H(u_n(t))\d t \\
  &\overset{\labelcref{eq6}}{=} -\intlim{G^-_n} \underbrace{H(u_n(t))}_{\geq -H^-}\d t\\
 &\overset{\hphantom{\labelcref{eq6}}}{\leq} H^- \underbrace{\mu(G^-_n)}_{\leq 1}.
\end{align*}
This implies
\begin{equation}\label{eq4}
\mu(\widetilde G_n) \leq \frac{\eps H^-}{R_n-1}.
\end{equation}
Now, let $t_n'\in \widetilde G_n$ be so that $[t_n',t_n]\subset\widetilde G$ and $r_n(s_n,t_n')=R_n-1$. Then the mean value theorem implies the existence of some $\theta_n \in [t_n',t_n]$ with \begin{align*}
\partial_t r_n(s_n,\theta_n) &= \frac{r_n(s_n,t_n)-r_n(s_n,t_n')}
{t_n-t_n'}\\
&= \frac{R_n-(R_n-1)}{\mu([t_n',t_n])}\\
&\geq \frac{1}{\mu(\widetilde G_n)}\\
&\overset{\labelcref{eq4}}{\geq} \frac{R_n-1}{\eps H^-}.
\end{align*} 
The desired inequality for $\partial_t u_n(s_n,\theta_n)$ then follows from Lemma \ref{lem3}.\qed\medskip

Assuming \labelcref{eqBoundConstraint}, we can now show that the cylinders in $\M(x^-,x^+)$ satisfy the conditions of \aref{lem3.4}:

\lemma\label{thm3.3} There is a compact $K\subset M$ with $\im(u)\subset K$ for all $u\in\M(x^-,x^+)$.
\proof We assume the contrary: Let $\folge{u}{n}\subset \M(x^-,x^+)$ be unbounded, i.e. let there be a sequence $(s_n,t_n)_{n\in\N}\subset\R\times\S$ and $R_n \vcentcolon = r_n(s_n,t_n)$ with $\lim\limits_{n\nach\infty} R_n=\infty$ as in \aref{lem3.3}.

First, a rescaling by $\frac{1}{\alpha_n}$ (with $\alpha_n \vcentcolon = \frac{R_n-1}{H^-}$ as in \aref{lem3.3}) for all $n\in\N$ yields
\begin{equation*}
\tilde u_n(s,t) \vcentcolon = u_n\left(s_n+\frac{s}{\alpha_n},t_n+\frac{t}{\alpha_n}\right).
\end{equation*}
Then, after writing $\tilde u_n=(\tilde r_n,\tilde x_n)$, $C_n \vcentcolon = R_n-1$, and $z\vcentcolon =s+it\in\C$ instead of $(s,t)\in\rs$, we set
\begin{equation*}
 w_n(z) \vcentcolon = \begin{cases} (\tilde r_n(z) - C_n, \tilde x_n(z)) &: \tilde r_n(z)-C_n \geq 0,\\
(0,\tilde x_n(z)) &: \text{  else},\end{cases}
\end{equation*}
so that $r(w_n(0)) = 1$ and $r(w_n(z)) \in [0,1]$ for all $z\in\C$.
 
For $z\in U_n \vcentcolon = \{z\in\C: r(\tilde u_n(z))>C_n\}$, these functions are still smooth and further satisfy (the shift by $-C_n$ vanishes when differentiating)
\begin{align} \label{eqX}
\nonumber \partial_s w_n + J(w_n)(\partial_t w_n) &= \partial_s \tilde u_n +J(  w_n)\partial_t \tilde u_n\\
&=  -\frac{1}{\alpha_n}\mult(\tilde u_n)J(\tilde u_n)X_H(\tilde u_n)\\
&= -\frac{1}{\alpha_n}\mult(\tilde u_n)J(w_n)X_H(w_n)
\end{align}
since $X_H$ and the SFT-like almost complex structure $J$ are independent of $r$ in $[0,\infty)\times\Sigma$, i.e. $X_H(w_n) = X_H(\tilde u_n)$ and $J(w_n) = J(\tilde u_n)$.

Further, the $\folge{w}{n}$ still have uniformly bounded energy, and because their images are now contained in the compact set $[0,1]\times\Sigma\subset M$, \aref{lem3.4} ensures that their derivatives are also uniformly bounded, i.e. the sequence is \emph{equicontinuous}.

The right hand side  of \labelcref{eqX} is not only bounded due to \aref{lemBoundConstraint}  and \aref{corGrad}, but even tends to zero in the $C^0$-norm due to the choice of $(\alpha_n)_{n\in\N}$.

Thus, if a smooth limit curve $w$ (for some subsequence) existed, it would satisfy
\begin{equation*}
\partial_s w + J(w)\partial_t w = 0,
\end{equation*}
i.e. be a $J$-holomorphic curve. To show the existence of $w$, one first uses the $L^\infty$ boundedness and equicontinuity of $\folge{w}{n}$ to apply Arzela-Ascoli, which gives weak convergence in $W^{1,\infty}$. Then, in order to show smoothness, \labelcref{eqX} is used for an elliptic bootstrapping argument. For details, see e.g. \cite{McDuffSalamon}, Theorem B.4.2, p. 575.

The only question remaining is  whether or not the domain of $w$ is non-empty. Here, we again use equicontinuity: Let $\delta>0$ be a bound on the (radial) derivatives of the $\folge{w}{n}$, then the following is true:
\newpage
\begin{claim}  $U\vcentcolon =B_{\frac{1}{\delta}}(0)\subset U_n$ for all $n\in\N$.\end{claim} 
\proof  
Let $z\in \C$ be arbitrary, then the mean value theorem and $r(\tilde u_n(0))=R_n$ imply the existence of a $\lambda\in [0,1]$ with \begin{align*} R_n-r(\tilde u_n(z)) &\leq \underbrace{\norm{\nabla r(\tilde u_n(\lambda z))}}_{\leq \delta}\cdot\norm{z}\leq \delta\norm{z}\end{align*}or equivalently $r(\tilde u_n(z)) \geq R_n-\delta\norm{z}$. For $\norm{z}\leq \frac{1}{\delta}$, this implies $z\in U$ which proves the claim for $\tilde u_n$. But since this inequality is invariant under the global shift by $-C_n$, the same property also holds for $w_n$.\qed

Now we apply \aref{lem3.3} with $\eps = \frac{1}{2\delta}$ to $\folge{u}{n}$. After rescaling by $\frac{1}{\alpha_n}$ and reparametrizing, the second inequality becomes \begin{equation}\label{eq22} \partial_t r(w_n(z_n))= \partial_t \tilde u_n(z_n) \geq \frac{1}{\eps}\end{equation}
for $z_n \vcentcolon = 0+i\alpha_n(\theta_n-0)$ and a sequence $\folge{\theta}{n}\subset\S$ with $$\alpha_n |\theta_n-0| \leq \frac{1}{2\delta}.$$

This means that $\folge{z}{n}$ is contained in a compact subset of $U$ for $n$ sufficiently large and thus converges to some $z_0\in U$. From \labelcref{eq22}, it then follows that $\partial_t r(w(z_0))\neq 0$, which implies that $w$ is non-constant in $U$. However, $w$ attains  a local maximum at $z=0$, and a computation as in in \cite{CFO}, Lemma 4.1, p. 32 shows that $r(w)$ is subharmonic, so the strong maximum principle leads to a contradiction. \qed
 
\subsection{Compactness up to Breaking}
As usual within Floer homology theories, one cannot expect a sequence of gradient flow cylinders to converge to another such cylinder: To compactify $\M(x^-,x^+)$, we need to add the following objects:

\definition A \emph{k-fold broken gradient flow cylinder} $(u^0,...,u^k)$   consists of $(k+1)$ gradient flow cylinders $u^i \in\M(x^-_i,x^+_i)$ with $x^+_i = x^-_{i+1}$ for all $i\in\{0,...,k-1\}$. If $k=0$, this reduces to a regular gradient flow cylinder.\\

With this terminology, we can state the final result of this section:

\theorem[Compactness]\label{thmLocGromov} If \labelcref{eqBoundConstraint} holds, then for every $\folge{u}{n}\subset \M(x^-,x^+)$ there exists a subsequence (denoted again by $\folge{u}{n}$) and a $k$-fold broken gradient flow cylinder $u=(u^0,...,u^k)$ with $x^-=x^-_0$ and $x^+_k=x^+$, so that (up to a reparametrization in the $s$-variable) $$u_n \overset{C^\infty_{loc}}{\nach} u^i$$ for some $i\in\{0,...,k\}$. Further, $k\leq \dim\M(x^-,x^+)$.
 
\proof This is a standard result for bounded lower order perturbations of $J$-holomorphic curves (cf. e.g. \cite{McDuffSalamon}, Theorem B.4.2, p. 575), and is proved by using the $L^\infty$-bounds on $u_n$, $\partial_s u_n$, $\partial_t u_n$, and $\mult(u_n)$ that were established in \aref{thm3.3}, \aref{lem3.4} and \aref{lemBoundConstraint}.\qed

\begin{remark} The result from \aref{thmLocGromov} is called \emph{compactness up to breaking}.\end{remark}

\newpage

\section{Definition of CFH}\label{secHom}
Before defining the new Floer homology, let us recall all assumptions used in this paper. See also the remarks after \aref{lemBoundConstraint} and \aref{corMF} for short discussions.\vspace{-1 \baselineskip}
\begin{assumption}\label{assC} In addition to \aref{assSymp}, let $H\in C^\infty(M,\R)$ be so that \labelcref{eqBoundConstraint} holds and assume that   surjectivity  in \aref{corMF} is satisfied.  
\end{assumption}
As in Rabinowitz Floer homology, the Morse-Bott situation necessitates the choice of an additional Morse function $h\in C^\infty(\crit{\act{H}},\R)$ and the following definition:
\begin{definition}  The Floer chain group $CF(\act{H},h)$ is the $\Z_2$-vector space of formal sums\begin{equation*}
\sum\limits_{x \in \crit{h}} \alpha_x x
\end{equation*} satisfying the \emph{Novikov semi-finiteness condition}
$$\bet{\menge{x \in\crit{h} :   \act{H}(x)\geq C\text{ and }\alpha_x\neq 0}} < \infty$$ for all $C\in\R$.
\end{definition}

Next, we need the concept of \emph{gradient flow lines with cascades}, which are explained in great detail in Appendix A of \cite{ExactContact}: Essentially, one chooses an auxiliary Riemannian metric $g$ on $\crit{\act{H}}$ so that $(h,g)$ is a Morse-Smale pair, and then considers flow lines from $x^-\in\crit{h}$ to $x^+\in\crit{h}$: If $x^\pm$ lie in the same connected component of $\crit{\act{H}}$, these are just the usual negative gradient flow lines of $h$ used to define Morse homology. Else, flow lines on different components are connected by negative gradient flow cylinders of $\act{H}$ called \emph{cascades}.

\definition Let $x^\pm\in\crit{h}$, and let $\Mb(x^-,x^+)$ be the moduli space of gradient flow lines with cascades between them. Then, the \emph{boundary operator} $\partial:CF(\act{H},h)\nach CF(\act{H},h)$ is the $\Z_2$-linear map induced by
    $$\partial x \vcentcolon= \sum\limits_{y\in \crit{h}} \cnt{\Mb^0(x,y)}y,$$ where $\cnt{\cdot}$ denotes the $\mod 2$ count of a finite set and $\Mb^0(x,y)$ is the $0$-dimensional component of $\Mb(x,y)$. 

\begin{remark}
 Finiteness of $\Mb^0(x,y)$  follows from combining \aref{thmLocGromov} with the same result for  Morse trajectories. By also making use of \emph{gluing}, the \glqq usual\grqq{} arguments from Floer theory show that $\partial^2=0$, leading to the following: 
\end{remark}

\definition\label{defCFH} The \emph{mean value constrained Floer homology} is given by:
\begin{equation*}
    CFH (M,\Sigma,H,J;h,g)\vcentcolon  = \frac{\ker \partial }{\im \partial }
\end{equation*}
 
As one would expect, this Floer homology is independent of $H$ in the following sense: 
\theorem \label{thmInv} The homology in \aref{defCFH} is independent of $J$ and the Morse-Smale pair $(h,g)$. Further, if $(H_s)_{s\in [0,1]}$ is a smooth family of defining Hamiltonians for $\Sigma_s = H_s\inv(0)$ as in \aref{assC}, then $CFH(\act{H_0})\cong CFH(\act{H_1})$. We thus also write $CFH(M,\Sigma)$.

\proof Due to the presence of the constraint factor $\mult$ in the gradient flow equation, the standard continuation principle from Hamiltonian Floer homology has to be adapted as in \cite{ExactContact}, p. 275 $-$ 277: The essential ingredient is an $L^\infty$-bound on $\mult$ which only depends on the asymptotics. In our case, this is \aref{lemBoundConstraint}. \qed

Using the grading   given by the (shifted)\footnote{For the definition of  homology, a global index shift is irrelevant since it cancels out in all formulas. When it comes to higher algebraic structures however, this shift  becomes more interesting. The next paper will contain a short discussion. } \emph{Morse-Bott index} $\bott $ introduced in \cite{ExactContact}, we can turn the $CF(\act{H},h)$ into a chain complex and refine above definition to 
\begin{equation*}
    CFH_* (M,\Sigma)\vcentcolon  = \frac{\ker \partial_* }{\im \partial_{*+1} }
\end{equation*}

\begin{remark}
Following our discussion in \textup{\cref{secRab}}, it is obvious that $$CF (\act{H})= CF (\Act{H}).$$ The significantly harder question $$CFH(M,\Sigma) \overset{?}{=} RFH(M,\Sigma)$$ will be explored in a future paper, using an adiabatic limit construction similar to that of Frauenfelder and Weber in \textup{\cite{Adiabatic}}.
\end{remark}

\newpage

\appendix

\section{Auxiliary Proofs}\label{secAuxProofs}
For a better reading flow, some proofs have been moved to this section.\vspace{\baselineskip}

\textit{Proof of }\aref{lemDiff}. Let $\gamma\in\loops[p]{M}$. It can be assumed to be continuous, so $t\auf H(\gamma(t))$ is continuous too and $$h(\gamma) = \int\limits_0^1 H(\gamma(t))\d t$$ is well-defined. For $\act{}$, we use that $\lambda = \omega(\liou,\cdot) = g(J\liou,\cdot)$ and \begin{align}\nonumber\act{}(\gamma) &= -\int\limits_\S \lambda_{\gamma(t)}(\partial_t\gamma) \\
&=-\int\limits_\S g_{\gamma(t)}(J\liou(\gamma(t)),\partial_t\gamma(t)).\label{eqA1}\end{align}
So even though $t\auf\partial_t\gamma(t)$ may not be well-defined, the integral \labelcref{eqA1} naturally is.

For their differentials, let $v\in W^{1,p}((-\eps_0,\eps_0),\loops{M}),\ \eps\auf v_\eps$ with $v_0=\gamma$ be a (weak) variation of $\gamma$. By construction of $\loops{M}$ via charts, there exists a unique \emph{continuous} vector field $X\in W^{1,p}(\S,\gamma\star\tang{}{M})$ with $v_\eps = \exp_\gamma \eps X$, i.e. $$\dbei{\eps}{0} v_\eps =X$$ in the \emph{classical} sense. Thus we can say:

\begin{align*}
\dbei{\eps}{0} h(v_\eps)&=\dbei{\eps}{0}\int\limits_0^1 H(v_\eps(t))\d t\\ 
&= \int\limits_0^1 \pbei{\eps}{0}H(v_\eps(t))\d t\\
&= \int\limits_0^1 H\left(\pbei{\eps}{0}v_\eps(t)\right)\d t\\
    &= \int\limits_0^1 \d H_{\gamma(t)} (X(t))\d t
\end{align*}

The calculation for $\dbei{\eps}{0}\act{}(v_\eps)$ works the same but involves the Lie derivative and Cartan's magic. Again, the \emph{classical} differentiability of $\eps\auf v_\eps$ is essential.\qed

\newpage 
The following should be read with the setting of \aref{thmSub} in mind, but it holds  for general Banach manifolds:
\theorem\label{thmBanachSub} Let $\L$ be a Banach manifold, $h\in C^1(\L,\R)$, and assume that at every $\gamma\in\H\vcentcolon =h\inv(0)$, there exists $v\in\tang{\gamma}{\L}$ with $\d h_\gamma(v) \neq 0$. Then $\H$ is a codimension $1$ submanifold of $\L$ and $\tang{\gamma}{\H}=\ker\d h_\gamma$. 
\proof
The assumption on $\d h$ already implies that every $\gamma\in h\inv(0)$ is a regular value. To use the infinite-dimensional version of the \emph{regular value theorem}\footnote{See e.g. \cite{Lang}, Chapter II.2, Proposition 2(ii).}, we still need to show that $\ker\d h$ \emph{splits} at every such $\gamma$, i.e.
$$\tang{\gamma}{\L} = \ker \d h_\gamma \oplus V_\gamma$$ for some closed (or just finite-dimensional) subspace $V_\gamma$. To this end, let $v\in \tang{\gamma}{\L}$ with $\d h_\gamma(v)\neq 0$, then
$$\pi : \tang{\gamma}{\L}\nach\tang{\gamma}{\L},\ w\auf w-  \frac{\d h_\gamma(w)}{\d h_\gamma(v)}v$$
is a linear map satisfying
\begin{enumerate}
    \item $\ker\pi = \langle v\rangle$,
    \item $\im\pi = \ker\d h_\gamma$,
    \item $\pi |_{\ker\d h_\gamma} = \id|_{\ker\d h_\gamma}$,
\end{enumerate} in other words $$\tang{\gamma}{\L} = \ker\d h_\gamma\oplus \langle v\rangle$$ via $$w = \pi(w) + v\frac{\d h_\gamma(w)}{\d h_\gamma(v)}$$ for $w\in \tang{\gamma}{\L}$. \qed

\remark In the proof of \aref{thmBanachSub}, every choice of $v\in\tang{\gamma}{\L}$ with $\d h_\gamma(v)\neq 0$ amounts to a different projection $\pi_v : \tang{\gamma}{\L}\nach\tang{\gamma}{\H}$. However, a right-inverse to all such projections is given by the identification map $$\iota : \tang{\gamma}{\H}\nach\tang{\gamma}{\L},\ w\auf w.$$

\theorem\label{thmGradBanach} Let $(\L,g)$ be a Hilbert manifold, $h\in C^1(\L,\R)$ have $0$ as a regular value, and $\H\vcentcolon=h\inv(0)$. Then $\tilde g \vcentcolon= g(\iota\cdot,\iota\cdot)$ defines a canonical scalar product on $\H$, and for any $f\in C^\infty(\L,\R)$ with $\tilde f \vcentcolon= f|_\H$, one has
\begin{equation*}\nabla \tilde f = \nabla f + \mult\nabla h,\end{equation*}
where $\mult : \H\nach\R$ is given by $$\mult(\gamma) \vcentcolon = -\frac{\d f_\gamma(\nabla h (\gamma))}{\d h_\gamma(\nabla h(\gamma))}.$$

\proof The first part is trivial. For the second part, we note that on $\tang{\gamma}{\H}$, $$\d \tilde f_\gamma = \d f_\gamma \circ\iota,$$ and so
\begin{align*}
    \tilde g(\nabla \tilde f,\cdot) &= \d \tilde f \\
    &= \d f_\gamma \circ\iota \\
    &= g(\nabla f, \iota\cdot) \\
    &= g(\iota\pi_v(\nabla f)+ (\nabla f-\iota \pi_v(\nabla f)),\iota \cdot)\\
    &= \tilde g(\pi_v(\nabla f),\cdot) + \frac{\d h (\nabla f)}{\d h (v)}g(v,\iota\cdot).
\end{align*} The second term vanishes, as soon as we choose a $v\in\tang{\gamma}{\L}$ with $\d h_\gamma(v)\neq0$ that is $g$-orthogonal to $\tang{\gamma}{\H}=\ker \d h$; for example $v = \nabla h(\gamma)$. In other words, \begin{align*}\nabla \tilde f &= \pi_{\nabla h}(\nabla f)\\
&= \nabla f - \frac{\d h(\nabla f)}{\d h(\nabla h)}\nabla h\\
&= \nabla f - \frac{\d f(\nabla h)}{\d h(\nabla h)}\nabla h.\end{align*}\qed

\newpage

\printbibliography
\end{document}